\input ./style/arxiv-ba.cfg
\documentclass[ba,linksfromyear,preprint]{imsart}
\makeatletter
   \@ifpackageloaded{natbib}{}{\usepackage{natbib}}
\makeatother

\pubyear{2015}
\volume{10}
\issue{3}
\firstpage{737}
\lastpage{739}
\doi{10.1214/15-BA944A}

\newcommand{\enquote}[1]{``#1''}

\begin{document}

\begin{frontmatter}
\title{Comment on Article by Ferreira and Gamerman\thanksref{T1}}
\runtitle{Comment on Article by Ferreira and Gamerman}

\relateddois{T1}{Main article DOI: \relateddoi[ms=BA944]{Related item:}{10.1214/15-BA944}.}

\begin{aug}

\author[addr1]{\fnms{Michael} \snm{Chipeta}\thanksref{t2}\ead[label=e1]{Michael.Chipeta@liverpool.ac.uk}}%
\thankstext{t2}{Michael Chipeta is supported by the UK Economic and Social Research Council
through the award of a PhD Studentship.}
\and
\author[addr2]{\fnms{Peter J.} \snm{Diggle}\corref{}\ead[label=e2]{p.diggle@lancaster.ac.uk}}

\runauthor{M. Chipeta and P. J. Diggle}

\address[addr1]{Institute of Infection and Global Health, University of Liverpool, \printead{e1}}
\address[addr2]{Institute of Infection and Global Health, University of Liverpool, \printead{e2}}

\end{aug}

\end{frontmatter}


\section{Introduction}

This paper is a welcome addition to the growing literature on preferential
sampling in a geostatistical setting. Earlier papers cited by the authors have
shown that preferential sampling materially affects parameter estimation and
prediction. The authors now demonstrate that the same applies to design, or
more specifically to the optimal augmentation of an initial set of geostatistical
data that has been sampled preferentially. Almost in passing, the paper also sets out
an algorithm for Bayesian inference under preferential sampling that is a useful contribution
in its own right. Might we look forward to an {\tt R} package implementation of this?

Our comments fall into two categories:
theoretical
remarks on what we call {\it adaptive design}, including an explanation
of why this does not necessarily require
consideration of preferential sampling issues; practical
constraints that may limit the scope for theoretically optimal designs to
be used in practice, especially in low-resource settings.\vspace*{-2pt}

\section{Adaptive geostatistical design and preferential sampling}

The topic of geostatistical design is multi-faceted. One useful distinction is between adaptive and non-adaptive designs. A {\it non-adaptive} design is one that is completely determined
before any data are collected. An {\it adaptive} design is one in which an initial design
is augmented in a way that depends on the analysis of interim data.  We make two theoretical
comments that follow from the definition of preferential sampling
given in Diggle, Menezes and Su (\citeyear{DMS}).

Firstly, an adaptive design need not be preferential.
To see why, it is sufficient to consider a two-stage adaptive design, $X = (X_0,X_1)$ with associated measurement data $(Y_0,Y_1)$,
where subscripts 0 and 1 identify initial and follow-up stages, respectively. Similarly,
write $S=(S_0,S_1)$ for the corresponding decomposition of the latent process~$S$. Quite generally,
we can factorise the joint distribution of $(X,Y,S)$ as
\begin{eqnarray}\label{eq:decomposition}
[X,Y,S] &=& [S,X_0,Y_0,X_1,Y_1]\nonumber \\
&=& [S]  [X_0|S]  [Y_0|X_0,S]  [X_1|Y_0,X_0,S] [Y_1|X_1,Y_0,X_0,S].
\end{eqnarray}
On the right-hand side of (\ref{eq:decomposition}), if the initial design is
non-preferential,
$[X_0|S] = [X_0]$, whilst by construction $[X_1|Y_0,X_0,S] = [X_1|Y_0,X_0]$. It then follows
that
\begin{eqnarray}
[X,Y] & = & [X_0][X_1|X_0,Y_0] \times \int_S [Y_0|X_0,S] [Y_1|X_1,Y_0,X_0,S] [S] dS \nonumber \\
      & = & [X|Y_0] \times [Y|X]
\label{eq:product}
\end{eqnarray}
and the log-likelihood is a sum of two components, $\log [X|Y_0] + \log [Y|X]$.
This shows
that the conditional likelihood, $[Y|X]$, can legitimately be used for inference although,
depending on how $[X|Y_0]$ is specified, to do so may be inefficient.  The argument leading to (\ref{eq:product}) is closely related to the proof that if data are ``missing at random'' the
missingness mechanism can be ignored when using likelihood-based inference (Rubin, \citeyear{Rubin}), and
extends to multi-stage adaptive designs with essentially only notational changes.

Secondly, shared dependence of a design $X$ and the latent process $S$
on observed covariates does not necessarily render $X$ preferential. Specifically,
if $Z$ denotes the covariate process, then $[X,S|Z]= [S|Z][X|S,Z]$. The requirement for
the design to be
non-preferential is that $[X|S,Z]=[X|Z]$, which in general is
a  weaker requirement than $[X|S]=[X]$. This illustrates, not for the first time, that
spatial
statistical inference can be greatly simplified by judicious selection of spatially referenced
covariates.

\section{Some practical constraints on geostatistical design}

The paper makes a number of explicit and
implicit assumptions that together provide a very reasonable
framework for theoretical analysis, but it is worth bearing in mind that
in any particular application, the design problem may be constrained in various ways.
These assumptions include the following:
\begin{enumerate}
   \item {\it The spatial integral of the predictive variance is an appropriate measure of
   predictive performance}

   This would not be true if, for example, $S(x)$ represents pollution and the main objective is to
   monitor compliance with environmental standards; see Fanshawe and Diggle (\citeyear{FD}).

  \item {\it Sampling may not be equally costly at every location}

  Put another way, should the design be constrained by the number of locations to be sampled,
or by the total sampling effort in the field?  An obvious example of this
is when travel-time represents a non-negligible proportion of field-effort; see,
for example, Figures 2 and 4
of \citet{DTetc}.

 \item
   {\it The number of potential sampling points may be finite}

   This applies to disease prevalence surveys when the sampling unit is either a household or a
   well-defined community. We are currently working on the adaptive design of
   an ongoing malaria prevalence mapping project around the perimeter of the Majete
   national park,  Malawi, where the first task has been to
   enumerate and geo-locate each household in each village within the study-region.  In the
   course of the project, we expect to sample all households, but the order in which they are
   sampled (in a sequence of monthly field-trips)
   will be chosen adaptively with the aim of optimising the estimation of
   the complete spatio-temporal variation in malaria prevalence, which is known to include
   a strong seasonal component.
\end{enumerate}

None of these these comments are intended to detract from the value of the paper
on its own terms.
Theoretical studies of this kind
help to further our understanding of important, and often
subtle,  methodological issues around modelling and inference for preferentially
sampled geostatistical data.

\end{document}